\documentclass{amsart}
\usepackage{amssymb}
\def\version{Version W01-v1m	 last changed 30 March 2004 by PG}
\date{\version}
\usepackage{a4}
\begin{document}
\newtheorem{theorem}{Theorem}[section]
\newtheorem{lemma}[theorem]{Lemma}
\newtheorem{example}[theorem]{Example}
\def\ffrac#1#2{{\textstyle\frac{#1}{#2}}}
\def\qedbox{\hbox{$\rlap{$\sqcap$}\sqcup$}}
\def\Rank{\operatorname{Rank}}
\makeatletter
 \renewcommand{\theequation}{%
 \thesection.\alph{equation}}
 \@addtoreset{equation}{section}
 \makeatother
\title[Covariant derivative curvature tensors]
{The structure of algebraic covariant derivative curvature tensors }
\author{J. D\'{\i}az-Ramos, B. Fiedler, E. Garc\'{\i}a-R\'{\i}o, and P. Gilkey}
\begin{address}{JD: Department of Geometry and Topology, Faculty of Mathematics, University of Santiago de Compostla,
15782 Santiago de Compostela, Spain. Email: {\it xtjosec@usc.es}}\end{address}
\begin{address}{BF: Mathematics Institute, University of Leipzig, Augustusplatz 10/11,\newline
04109 Leipzig,
Germany. Email: {\it bernd.fiedler.roschstr.leipzig@t-online.de}}\end{address}
\begin{address}{EG: Department of Geometry and Topology, Faculty of Mathematics, University of Santiago de Compostla,
15782 Santiago de Compostela, Spain. Email: {\it xtedugr@usc.es}}\end{address}
\begin{address}{PG: Mathematics Department, University of Oregon,
Eugene Or 97403 USA.\newline Email: {\it gilkey@darkwing.uoregon.edu}}
\end{address}

\begin{abstract} We use the Nash embedding theorem to construct generators for the space of
algebraic covariant derivative curvature tensors.\end{abstract}
\maketitle

\section{Introduction} Let $M$ be an $m$ dimensional Riemannian manifold. To a large extent, the geometry of $M$ is
the study of the Riemannian curvature $R\in\otimes^4T^*M$ which is defined by the Levi-Civita connection
$\nabla$ and, to a lesser extent, the study of the covariant derivative $\nabla R$. For example, $M$ is a
local symmetric space if and only if $\nabla R=0$; note that local symmetric spaces are locally homogeneous.

It is convenient to work in the algebraic context. Let $V$ be an $m$-dimensional real vector space.
Let $\mathcal{A}(V)\subset\otimes^4V^*$ and $\mathcal{A}_1(V)\subset\otimes^5V^*$ be the spaces of all algebraic curvature tensors
and all algebraic covariant derivative tensors, respectively, i.e. those tensors $A$ and $A_1$ having the symmetries of $R$ and of
$\nabla R$:
\begin{eqnarray*}
&&A(x,y,z,w)=A(z,w,x,y)=-A(y,x,z,w),\\
&&A(x,y,z,w)+A(y,z,x,w)+A(z,x,y,w)=0,\\
&&A_1(x,y,z,w;v)=A_1(z,w,x,y;v)=-A_1(y,x,z,w;v),\\
&&A_1(x,y,z,w;v)+A_1(y,z,x,w;v)+A_1(z,x,y,w;v)=0,\\
&&A_1(x,y,z,w;v)+A_1(x,y,w,v;z)+A_1(x,y,v,z;w)=0\,.
\end{eqnarray*}
Let $S^p(V)\subset\otimes^pV^*$ be the space of totally symmetric $p$ forms. If $\Psi\in S^2(V)$ and if  $\Psi_1\in S^3(V)$, define
$A_\Psi\in\mathcal{A}(V)$ and $A_{1,\Psi,\Psi_1}\in\mathcal{A}_1(V)$ by:
\begin{eqnarray*}
A_\Psi(x,y,z,w):&=&\Psi(x,w)\Psi(y,z)-\Psi(x,z)\Psi(y,w),\\
A_{1,\Psi,\Psi_1}(x,y,z,w;v):&=&\Psi_1(x,w,v)\Psi(y,z)+\Psi(x,w)\Psi_1(y,z,v)\\
  &-&\Psi_1(x,z,v)\Psi(y,w)-\Psi(x,z)\Psi_1(y,w,v)\,.
\end{eqnarray*}
If one thinks of $\Psi_1$ as the symmetrized covariant derivative of $\Psi$, then $A_{1,\Psi,\Psi_1}$ can be regarded, at least
formally speaking, as the covariant derivative of $A_\Psi$.

Fiedler \cite{F01,F03} used group representation theory
to show:

\begin{theorem}[Fiedler]\label{thm-1.1}\ \begin{enumerate}
\item $\mathcal{A}(V)=\operatorname{Span}_{\Psi\in S^2(V)}\{A_\Psi\}$.
\item $\mathcal{A}_1(V)=\operatorname{Span}_{\Psi\in S^2(V),\Psi_1\in S^3(V)}\{A_{1,\Psi,\Psi_1}\}$.
\end{enumerate}
\end{theorem}

Let $A\in\mathcal{A}(V)$ and $A_1\in\mathcal{A}_1(V)$ be given. Choose $\nu(A)$ and $\nu_1(A_1)$ minimal so that there exist
$\Psi_i\in S^2(V)$, $\tilde\Psi_j\in S^2(V)$, $\tilde\Psi_{1,j}\in S^3(V)$, and constants $\lambda_i,\lambda_{1,j}$ so:
$$A=\textstyle\sum_{1\le i\le\nu(A)}\lambda_iA_{\Psi_i}\quad\text{and}\quad
A_1=\textstyle\sum_{1\le j\le\nu_1(A_1)}\lambda_{1,j}A_{1,\tilde\Psi_j,\tilde\Psi_{1,j}}\,.
$$
Set 
$$\nu(m):=\sup_{A\in\mathcal{A}(V)}\nu(A)\quad\text{and}\quad\nu_1(m):=\sup_{A_1\in\mathcal{A}_1(V)}\nu_1(A_1)\,.$$
The main result of
this paper is the following:

\begin{theorem}\label{thm-1.2}
Let $m\ge2$. \begin{enumerate}
\item $\frac12m\le\nu(m)$ and $\frac12m\le\nu_1(m)$.
\item $\nu(m)\le\ffrac12m(m+1)$ and $\nu_1(m)\le\ffrac12m(m+1)$.
\end{enumerate}\end{theorem}

 We shall establish the lower bounds of Assertion (1) in Section \ref{sect-2}.
The upper bound given in Assertion (2) for $\nu(m)$ is due to  D\'{\i}az-Ramos and
Garc\'{\i}a-R\'{\i}o \cite{DRGR-04} who used the Nash embedding theorem \cite{Na56}; they also gave a separate argument to show
$\nu(2)=1$ and $\nu(3)=2$. In Section
\ref{sect-3}, we shall generalize their approach to establish the following simultaneous `diagonalization' result from which
Theorem \ref{thm-1.2} (2) will follow as a Corollary:

\begin{theorem}\label{thm-1.3} Let $V$ be an $m$ dimensional vector space. Let $A\in\mathcal{A}(V)$ and let
$A_1\in\mathcal{A}_1(V)$ be given. There exists $\Psi_i\in S^2(V)$ and $\Psi_{1,i}\in S^3(V)$ so
that
$$A=\textstyle\sum_{1\le i\le\frac12m(m+1)}A_{\Psi_i}\quad\text{and}\quad
  A_1=\textstyle\sum_{1\le i\le\frac12m(m+1)}A_{1,\Psi_i,\Psi_{1,i}}\,.
$$
\end{theorem}

The study of the tensors $A_\Psi$ arose in the original instance from the Osserman conjecture and related matters; we refer to
\cite{GKV02,Gi02} for a more extensive discussion than is possible here, and content ourselves with only a very brief introduction
to the subject.

\subsection{The Jacobi operator} If $M$ is a pseudo-Riemannian manifold of signature $(p,q)$ and dimension $m=p+q$, let $S^+(M)$
(resp.
$S^-(M)$) be the bundle of unit spacelike (resp. timelike) tangent vectors. The Jacobi operator $J(x)$ for $x\in TM$ is the
self-adjoint endomorphism of $TM$ characterized by the identity:
$$g(J(x)y,z)=R(y,x,x,z)\,.$$
One says that $M$ is {\it spacelike Osserman} (resp. {\it timelike Osserman})
if the eigenvalues of $J(\cdot)$ are constant on $S^+(M)$ (resp. $S^-(M)$). It turns out these two notions are equivalent
and such a manifold is simply said to be {\it Osserman}. 

Restrict for the moment to the Riemannian setting ($p=0$). If $M$ is a local rank
$1$ symmetric space or is flat, then the local isometries of $M$ act transitively on the sphere bundle $S(M)=S^+(M)$ and
hence the eigenvalues of $J(\cdot)$ are constant on $S(M)$ and $M$ is Osserman. Osserman \cite{O97}
wondered if the converse held; this question has been called the Osserman conjecture by subsequent authors. The conjecture has been
answered in the affirmative if
$m\ne16$ by work of Chi \cite{Ch88} and Nikolayevsky \cite{Ni03,Ni04a,Ni04b}. 

In the Lorentzian setting ($p=1$), an Osserman
manifold has constant sectional curvature
\cite{BBG97,GKV97}. In the higher signature setting ($p>1$, $q>1$) it is more natural to work with the Jordan
normal form rather than just the eigenvalue structure. One says that
$M$ is {\it spacelike Jordan Osserman} (resp. {\it timelike Jordan Osserman}) if the Jordan normal form of $J(\cdot)$ is
constant on
$S^+(M)$ (resp.
$S^-(M)$); these two notions are not equivalent. The following example is
instructive. Let
$(\vec x,\vec y)$ for $\vec x=(x_1,...,x_p)$ and $\vec y=(y_1,...,y_p)$ be coordinates on $\mathbb{R}^{2p}$ where $p\ge3$. Let
$f=f(\vec x)\in C^\infty(\mathbb{R}^p)$. Define a pseudo-Riemannian metric $g_f$ of signature $(p,p)$ on $\mathbb{R}^{2p}$ by
setting
\begin{equation}\label{eqn-1.c}
g_f(\partial_i^x,\partial_j^x)=\partial_i^xf\cdot\partial_j^xf,\ \ 
  g_f(\partial_i^y,\partial_j^y)=0,\ \text{ and }\ 
g_f(\partial_i^x,\partial_j^y)=g_f(\partial_j^y,\partial_i^x)=\delta_{ij}\,.
\end{equation}
Let $\Psi$ be the Euclidean Hessian:
$$\Psi(\partial_i^x,\partial_j^x)=\partial_i^x\partial_j^xf,\quad
\Psi(\partial_i^y,\partial_j^y)=0,\quad\text{and}\quad
\Psi(\partial_i^x,\partial_j^y)=\Psi(\partial_j^y,\partial_i^x)=0\,.
$$
One then has that $R=A_\Psi$. We suppose that the restriction of $\Psi$ to
$\operatorname{Span}\{\partial_i^x,\partial_j^x\}$ is positive definite henceforth. Then
$M$ is a complete pseudo-Riemannian manifold which is spacelike and timelike Jordan Osserman. Similarly set
$$\Psi_1(\partial_i^x,\partial_j^x,\partial_k^x)=\partial_i^x\partial_j^x\partial_k^xf$$
and extend $\Psi_1$ to vanish if any
entry is
$\partial_\ell^y$. One has
$\nabla R=A_{1,\Psi,\Psi_1}$; thus if $f$ is not quadratic, $M$ is not a local symmetric
space. With a bit more work one can show that for generic such $f$, $M$ is curvature homogeneous but
not locally affine homogeneous. We refer to \cite{DG04,GS04} for further details. 

\subsection{The skew-symmetric curvature operator} Let $\{e_1,e_2\}$ be an orthonormal basis for an oriented spacelike (resp.
timelike)
$2$ plane
$\pi$. The skew-symmetric curvature operator $\mathcal{R}(\pi)$ is characterized by the identity
$$g(\mathcal{R}(\pi)y,z)=R(e_1,e_2,y,z)\,;$$
it is independent of the particular orthonormal basis chosen. One says that $M$ is {\it spacelike Ivanov-Petrova} (resp.
{\it timelike Ivanov-Petrova}) if the eigenvalues of $\mathcal{R}(\cdot)$ are constant on the Grassmannian of oriented spacelike
(resp. timelike) $2$-planes; these two notions are equivalent and such a manifold is simply said to be Ivanov-Petrova. The
notions {\it spacelike Jordan Ivanov-Petrova} and {\it timelike Jordan Ivanov-Petrova} are defined similarly and are not
equivalent. 

The Riemannian Ivanov-Petrova manifolds have been classified \cite{Gi99,GLS99,Ni04c}; they have also been classified in the
Lorentzian setting \cite{Z02} if $m\ge10$. For all these manifolds, the curvature tensors have the form
$R=A_\Psi$ where $\Psi$ is an idempotent isometry and $\mathcal{R}(\pi)$ always has rank $2$. Conversely, in the algebraic
setting, if
$R$ is a spacelike Jordan Ivanov-Petrova algebraic curvature tensor on a vector space of signature $(p,q)$ where $q\ge5$ and where
$\Rank\{\mathcal{R}(\cdot)\}=2$, then there exist $\lambda$ and $\Psi$ so that
$R=\lambda A_\Psi$. This once again motivates the study of these tensors. Unfortunately, the situation in the indefinite
setting is again quite different. There exist spacelike Ivanov-Petrova manifolds of signature
$(s,2s)$ where $\mathcal{R}(\pi)$ has rank $4$ and where the curvature tensor does not have the form $R=A_\Psi$. We refer to
\cite{GS04a} for further details.

\subsection{The Szab\'o operator} There is an analogous operator to the Jacobi operator which is defined by $\nabla R$. The Szab\'o
operator
$J_1(x)$ is the self-adjoint endomorphism of
$TM$ characterized by $g(J_1(x)y,z)=\nabla R(y,x,x,z;x)$. One says that $M$ is {\it spacelike Szab\'o} (resp. {\it
timelike Szab\'o}) if the eigenvalues of $J_1(\cdot)$ are constant on $S^+(M)$ (resp. $S^-(M)$); these notions are equivalent and
such a manifold is simply said to be Szab\'o. The notion {\it spacelike} (resp. {\it timelike}) {\it Jordan Szab\'o} is defined
similarly.

In his study of $2$ point symmetric spaces, Szab\'o \cite{S91} gave a very lovely topological argument showing that any
Riemannian Szab\'o manifold is necessarily a local symmetric space -- i.e. $\nabla R=0$. This result was subsequently extended
to the Lorentzian case \cite{GSt02}. In the higher signature setting, again the situation is unclear. The metric $g_f$ described in
Display (\ref{eqn-1.c}) defines a Szab\'o pseudo-Riemannian manifolds of signature $(p,p)$. 

Even in the
algebraic setting, there are no known non-zero elements $A_1\in\mathcal{A}(V)$ which are spacelike Jordan Szab\'o. It has
been shown \cite{GIS03} that if $A_1$ is a spacelike Jordan Szab\'o algebraic covariant derivative curvature tensor on a vector
space of signature $(p,q)$, where
$q\equiv1$ mod $2$ and $p<q$ or where $q\equiv2$ mod $4$ and $p<q-1$, then $A_1=0$. This algebraic result yields an elementary
proof of the geometrical fact that any pointwise totally isotropic pseudo-Riemannian manifold with such a signature
$(p,q)$ is locally symmetric. The general question of finding non-trivial spacelike Jordan Szab\'o covariant algebraic curvature
tensors, or conversely showing non exist, remains open.

The examples discussed above motivate consideration of the tensors $A_{1,\Psi,\Psi_1}$ and more generally of tensors which are
combinations of these. We hope that Theorems \ref{thm-1.2} and \ref{thm-1.3}, although of interest in their own right, will play a
central role in these investigations.

\section{A lower bound for $\nu(m)$ and for $\nu_1(m)$}\label{sect-2}

Let $V$ be an $m$ dimensional vector space,  let
$A\in\mathcal{A}(V)$, and let $A_1\in\mathcal{A}_1(V)$. Give $V$ a positive definite inner product
$\langle\cdot,\cdot\rangle$. The associated curvature operators are then defined by the identities: 
\begin{eqnarray*}
&&\langle \mathcal{R}_A(\xi_1,\xi_2)z,w\rangle=A(\xi_1,\xi_2,z,w),\quad\text{and}\\
&&\langle \mathcal{R}_{A_1}(\xi_1,\xi_2,\xi_3)z,w\rangle=A_1(\xi_1,\xi_2,z,w;\xi_3)\,.
\end{eqnarray*}
Theorem \ref{thm-1.2} (1) will follow from the following Lemma:
\begin{lemma}\label{lem-2.1}Let $V$ be a vector space of dimension $m=2\bar m$ or $m=2\bar m+1$. \begin{enumerate}
\item If $\Psi\in S^2(V)$ and if $\Psi_1\in S^3(V)$, then for any $\xi_1,\xi_2,\xi_3\in V$ one has:
$$
  \Rank\{\mathcal{R}_{A_\Psi}(\xi_1,\xi_2)\}\le2\quad\text{and}\quad
  \Rank\{\mathcal{R}_{A_{1,\Psi,\Psi_1}}(\xi_1,\xi_2,\xi_3)\}\le2\,.
$$
\item If $A\in\mathcal{A}(V)$ and $A_1\in\mathcal{A}_1(V)$, then for any $\xi_1,\xi_2,\xi_3\in V$ one has:
$$
  \Rank\{\mathcal{R}_A(\xi_1,\xi_2)\}\le2\nu(A)\quad\text{and}\quad
  \Rank\{\mathcal{R}_{A_1}(\xi_1,\xi_2,\xi_3)\}\le2\nu_1(A_1)\,.
$$
\item There exist $A\in\mathcal{A}(V)$, $A_1\in\mathcal{A}_1(V)$, and $\xi_1,\xi_2,\xi_3\in V$ so:
$$
  \Rank\{\mathcal{R}_A(\xi_1,\xi_2)\}=2\bar m\quad\text{and}\quad
  \Rank\{\mathcal{R}_{A_1}(\xi_1,\xi_2,\xi_1)\}=2\bar m\,.
$$
\end{enumerate}\end{lemma}

\medbreak\noindent{\it Proof.} If $\Psi\in
S^2(V)$ and
$\Psi_1\in S^3(V)$, let
$\psi$ and
$\psi_1(\cdot)$ be the associated self-adjoint endomorphisms characterized by the identities
$$\langle\psi x,y\rangle=\Psi(x,y)\quad\text{and}\quad \langle\psi_1(z)x,y\rangle=\Psi_1(x,y,z)\,.$$
Assertion (1) follows from the expression:
\begin{eqnarray*}
  \mathcal{R}_{A_\Psi}(\xi_1,\xi_2)y&=&\{\Psi(\xi_2,y)\psi\}\xi_1-\{\Psi(\xi_1,y)\psi\}\xi_2,\quad\text{and}\\
  \mathcal{R}_{A_{1,\Psi,\Psi_1}}(\xi_1,\xi_2,\xi_3)y&=&\{\Psi(\xi_2,y)\psi_1(\xi_3)+\Psi_1(\xi_2,y,\xi_3)\psi\}\xi_1\\
&-&\{\Psi(\xi_1,y)\psi_1(\xi_3)+\Psi_1(\xi_1,y,\xi_3)\psi\}\xi_2\,.
\end{eqnarray*}
Let $A_i:=A_{\Psi_i}$, $A_{1,j}:=A_{1,\tilde\Psi_j,\tilde\Psi_{1,j}}$, $\mathcal{R}_i:=\mathcal{R}_{A_i}$, and
$\mathcal{R}_{1,i}:=\mathcal{R}_{A_{1,i}}$. Set
$$\textstyle A=\sum_{1\le
i\le\nu(A)}A_i\quad\text{and}\quad
A_1=\sum_{1\le j\le\nu_1(A_1)}A_{1,j}\,.
$$
Assertion (2) follows from Assertion (1) as
\begin{eqnarray*}
\Rank\{\mathcal{R}_A(\cdot)\}
&=&\Rank\{\textstyle\sum_{1\le i\le\nu(A)}\mathcal{R}_i(\cdot)\}\\
&\le&\textstyle\sum_{1\le i\le\nu(A)}
\Rank\{\mathcal{R}_i(\cdot)\}\le2\nu(A),\\
\Rank\{\mathcal{R}_{A_1}(\cdot)\}
&=&\Rank\{\textstyle\sum_{1\le j\le\nu_1(A_1)}\mathcal{R}_{1,j}(\cdot)\}\\
&\le&\textstyle\sum_{1\le j\le\nu_1(A_1)}\Rank\{\mathcal{R}_{1,j}(\cdot)\}\le2\nu_1(A_1)\,.
\end{eqnarray*}

If $\dim(V)=2\bar m$, let $\{e_1,...,e_{\bar m},f_1,...,f_{\bar m}\}$ be an orthonormal basis for $V$; if $\dim(V)$ is odd, the
argument is similar and we simply extend $A$ and $A_1$ to be trivial on the additional basis vector. Define the non-zero components
of
$\Psi_i\in S^2(V)$ and $\Psi_{1,i}\in S^3(V)$ by:
\begin{eqnarray*}
&&\Psi_i(e_j,e_k)=\Psi_i(f_j,f_k)=\delta_{ij}\delta_{ik},\\
&&\Psi_{1,i}(e_j,e_k,e_l)=\Psi_{1,i}(f_j,f_k,f_l)=\delta_{ij}\delta_{ik}\delta_{il};
\end{eqnarray*}
$\Psi_i(\cdot,\cdot)$ and $\Psi_{1,i}(\cdot,\cdot,\cdot)$ vanish if both an `e' and an `f'
appear. Let
\begin{eqnarray*}
&&A_i:=A_{\Psi_i},\quad\mathcal{R}_i:=\mathcal{R}_{A_i},\quad
  A_{1,i}:=A_{1,\Psi_i,\Psi_{1,i}},\quad\mathcal{R}_{1,i}:=\mathcal{R}_{A_{1,i}},\\
&&A:=\textstyle\sum_{1\le i\le\bar m}A_i,\quad A_1:=\textstyle\sum_{1\le i\le\bar m}A_{1,i},\\
&&\xi_1:=e_1+...+e_{\bar m},\quad
\xi_2:=f_1+...+f_{\bar m},\quad\xi_3:=\xi_1+\xi_2\,.
\end{eqnarray*}
We may then complete the proof of Assertion (3) by computing:
\begin{eqnarray*}
&&\mathcal{R}_A(\xi_1,\xi_2)e_i=\mathcal{R}_i(e_i,f_i)e_i=-f_i,\\
&&\mathcal{R}_A(\xi_1,\xi_2)f_i=\mathcal{R}_i(e_i,f_i)f_i=e_i,\\
&&\mathcal{R}_{A_1}(\xi_1,\xi_2,\xi_3)e_i=\mathcal{R}_{1,i}(e_i,f_i,e_i+f_i)e_i=-2f_i\\
&&\mathcal{R}_{A_1}(\xi_1,\xi_2,\xi_3)f_i=\mathcal{R}_{1,i}(e_i,f_i,e_i+f_i)f_i=2e_i\,.\quad\qedbox
\end{eqnarray*}

\section{Geometric realizability}\label{sect-3} Henceforth,  let $\langle\cdot,\cdot\rangle$ be a non-singular innerproduct on an
$m$ dimensional vector space $V$, let $A\in\mathcal{A}(V)$ and let $A_1\in\mathcal{A}(V)$. 

Although the following is well-known, see for
example Belger and Kowalski \cite{BK04} where a more general result is established, we shall give the proof to keep the development
as self-contained as possible and to establish notation needed subsequently.
\begin{lemma}\label{lem-3.1}\ \begin{enumerate}
\item If $g$ is a pseudo-Riemannian metric on $\mathbb{R}^m$ with
$\partial_ig_{jk}(0)=0$, then:
\begin{enumerate}\smallbreak\item
$R_{ijkl}(0)=\textstyle\frac12\{\partial_i\partial_kg_{jl}+\partial_j\partial_l g_{ik}
    -\partial_i\partial_l g_{jk}-\partial_j\partial_kg_{il}
    \}(0)$.
\item $R_{ijkl;n}(0)=\textstyle\frac12\{\partial_i\partial_k\partial_ng_{jl}+\partial_j\partial_l\partial_n g_{ik}
    -\partial_i\partial_l\partial_n g_{jk}-\partial_j\partial_k\partial_ng_{il}
    \}(0)$.\end{enumerate}
\smallbreak\item
There exists the germ of a pseudo-Riemannian
metric $g$ on $(\mathbb{R}^m,0)$ and an isomorphism $\Xi$ from $T_0(\mathbb{R}^m)$ to $V$ so that
\begin{enumerate}\smallbreak\item
$\Xi^*\langle\cdot,\cdot\rangle=g|_{T_0(\mathbb{R}^m)}$.
\item$\Xi^*A=R_g|_{T_0(\mathbb{R}^m)}$.
\item $\Xi^*A_1=\nabla R_g|_{T_0(\mathbb{R}^m)}$.
\end{enumerate}
\end{enumerate}\end{lemma}

\medbreak\noindent{\it Proof.} Since the $1$ jets of the metric vanish at the origin, we have
\begin{eqnarray*}
&&\Gamma_{ijk}:=g(\nabla_{\partial_i}\partial_j,\partial_k)=\ffrac12(\partial_ig_{jk}+\partial_jg_{ik}-\partial_kg_{ij})
=O(|x|),\\
&&R_{ijkl}(0)=\{\partial_i\Gamma_{jkl}-\partial_j\Gamma_{ikl}\}(0),\quad\text{and}\quad
R_{ijkl;n}(0)=\{\partial_nR_{ijkl}\}(0)\,.
\end{eqnarray*}
Assertion (1) now follows; see, for example, \cite{Gi02}
[cf Lemma 1.11.1] for further details. To prove the second assertion, choose an orthonormal basis
$\{e_1,...,e_m\}$ for
$V$ so that $\langle e_i,e_j\rangle=\pm\delta_{ij}$; we use this orthonormal basis to identify
$V=\mathbb{R}^m$. Let $A_{ijkl}$ and $A_{1,ijkl;n}$ denote the components of $A$ and of $A_1$, respectively. Define
\begin{eqnarray*}
g_{ik}&=&\langle e_i,e_k\rangle-\textstyle\frac13\textstyle\sum_{jl}A_{ijlk}x_jx_l
     -\frac16\textstyle\sum_{jln}A_{1,ijlk;n}x_jx_lx_n\,.
\end{eqnarray*}
Clearly $g_{ik}=g_{ki}$. As $g|T_0\mathbb{R}^m=\langle\cdot,\cdot\rangle$, $g$ is non-degenerate on some neighborhood of $0$. Since
the
$1$ jets of the metric vanish at
$0$ we have by Assertion (1) that
\begin{eqnarray*}
&&R_{ijkl}(0)\\&=&
 \textstyle\frac16\{-A_{jikl}-A_{jkil}-A_{ijlk}-A_{iljk}
    +A_{jilk}+A_{jlik}+A_{ijkl}+A_{ikjl}\}\\
  &=&\textstyle\frac16\{4A_{ijkl}-2A_{iljk}-2A_{iklj}\}
 =A_{ijkl},\\
&&R_{ijkl;n}(0)\\
&=&\textstyle\frac1{12}\{-A_{jikl;n}-A_{jkil;n}-A_{jnkl;i}-A_{jknl;i}-A_{jinl;k}-A_{jnil;k}\\
&-&A_{ijlk;n}-A_{iljk;n}-A_{inlk;j}-A_{ilnk;j}-A_{ijnk;l}-A_{injk;l}\\
&+&A_{jilk;n}+A_{jlik;n}+A_{jnlk;i}+A_{jlnk;i}+A_{jink;l}+A_{jnik;l}\\
&+&A_{ijkl;n}+A_{ikjl;n}+A_{inkl;j}+A_{iknl;j}+A_{ijnl;k}+A_{injl;k}\}\\
&=&\textstyle\frac1{12}\{(4A_{ijkl;n}-2A_{jkil;n}+2A_{jlik;n})+(-2A_{jnkl;i}-2A_{inlk;j})\\
&+&(-2A_{jinl;k}-2A_{ijnk;l})+(-A_{ilnk;j}-A_{jnil;k})\\
&+&(-A_{injk;l}-A_{jknl;i})
+(A_{jlnk;i}+A_{injl;k})+(A_{jnik;l}+A_{iknl;j})\}\\
&=&\ffrac1{12}\{6A_{ijkl;n}+2A_{ijkl;n}+2A_{ijkl;n}+A_{ilkj;n}+A_{jkli;n}
-A_{jlki;n}-A_{iklj;n}\}\\
&=&\ffrac1{12}\{10A_{ijkl;n}+2A_{ilkj;n}+2A_{ikjl;n}\}
=\ffrac1{12}\{10A_{ijkl;n}-2A_{ijlk;n}\}=A_{ijkl;n}
\,.\quad\qedbox
\end{eqnarray*}

We suppose the inner product $\langle\cdot,\cdot\rangle$ is positive definite henceforth. We apply the Nash embedding theorem
\cite{Na56} to find an embedding $f:\mathbb{R}^m\rightarrow\mathbb{R}^{m+\kappa}$ realizing the metric $g$ constructed in Lemma
\ref{lem-3.1}. By writing the submanifold as a graph over its tangent plane, we can choose coordinates $(x,y)$ on
$\mathbb{R}^{m+\kappa}$ where $x=(x_1,...,x_m)$ and $y=(y_1,...,y_\kappa)$ so that
$$f(x)=(x,f_1(x),...,f_\kappa(x))\quad\text{where}\quad df_\nu(0)=0\quad\text{for}\quad1\le \nu\le \kappa\,.$$
Since $f_*(\partial_i^x)=(0,...,1,...,0,\partial_i^xf_1,....,\partial_i^xf_\kappa)$, we have
$$g_{ij}(x)=\delta_{ij}+\textstyle\sum_{1\le\sigma\le\kappa}\partial_i^xf_\sigma\cdot\partial_j^xf_\sigma\,.$$
Let $\Psi^\sigma_{ij}:=\partial_i^x\partial_j^xf_\sigma(0)$ and
$\Psi^\sigma_{ijk}:=\partial_i^x\partial_j^x\partial_k^xf_\sigma(0)$.  As $dg_{ij}(0)=0$, by Lemma \ref{lem-3.1}:
\begin{eqnarray*}
&&R_{ijkl}(0)\\
&=&\ffrac12\textstyle\sum_{1\le\sigma\le\kappa}\{
(\Psi^\sigma_{ij}\Psi^\sigma_{kl}+\Psi^\sigma_{il}\Psi^\sigma_{kj})+(\Psi^\sigma_{ji}\Psi^\sigma_{lk}+\Psi^\sigma_{jk}\Psi^\sigma_{li})\\
&-&(\Psi^\sigma_{ij}\Psi^\sigma_{lk}+\Psi^\sigma_{ik}\Psi^\sigma_{lj})-(\Psi^\sigma_{ji}\Psi^\sigma_{kl}+\Psi^\sigma_{jl}\Psi^\sigma_{ki})\}\\
&=&\textstyle\sum_{1\le\sigma\le\kappa}\{\Psi^\sigma_{il}\Psi^\sigma_{jk}-\Psi^\sigma_{ik}\Psi^\sigma_{jl}\}
=\sum_{1\le\sigma\le\kappa} A_{\Psi^\sigma},\\
&&R_{ijkl;n}(0)\\
&=&\ffrac12\textstyle\sum_{1\le\sigma\le\kappa}\{
   (\Psi^\sigma_{jin}\Psi^\sigma_{lk}+\Psi^\sigma_{jkn}\Psi^\sigma_{li}+\Psi^\sigma_{ji}\Psi^\sigma_{lkn}+\Psi^\sigma_{jk}\Psi^\sigma_{lin}
          +\Psi^\sigma_{jik}\Psi^\sigma_{ln}+\Psi^\sigma_{jn}\Psi^\sigma_{lik})\\
&+&(\Psi^\sigma_{ijn}\Psi^\sigma_{kl}+\Psi^\sigma_{iln}\Psi^\sigma_{kj}+\Psi^\sigma_{ij}\Psi^\sigma_{kln}+\Psi^\sigma_{il}\Psi^\sigma_{kjn}
          +\Psi^\sigma_{ijl}\Psi^\sigma_{kn}+\Psi^\sigma_{in}\Psi^\sigma_{kjl})\\
&-&(\Psi^\sigma_{jin}\Psi^\sigma_{kl}+\Psi^\sigma_{jln}\Psi^\sigma_{ki}+\Psi^\sigma_{ji}\Psi^\sigma_{kln}+\Psi^\sigma_{jl}\Psi^\sigma_{kin}
          +\Psi^\sigma_{jil}\Psi^\sigma_{kn}+\Psi^\sigma_{jn}\Psi^\sigma_{kil})\\
&-&(\Psi^\sigma_{ijn}\Psi^\sigma_{lk}+\Psi^\sigma_{ikn}\Psi^\sigma_{lj}+\Psi^\sigma_{ij}\Psi^\sigma_{lkn}+\Psi^\sigma_{ik}\Psi^\sigma_{ljn}
          +\Psi^\sigma_{ijk}\Psi^\sigma_{ln}+\Psi^\sigma_{in}\Psi^\sigma_{ljk})\\
&=&\textstyle\sum_{1\le\sigma\le\kappa}\{\Psi^\sigma_{iln}\Psi^\sigma_{jk}+\Psi^\sigma_{jkn}\Psi^\sigma_{il}
-\Psi^\sigma_{ikn}\Psi^\sigma_{jl}-\Psi^\sigma_{ik}\Psi^\sigma_{jln}\}
=\textstyle\sum_{1\le\sigma\le\kappa} A_{1,\Psi^\sigma,\Psi^\sigma}\,.
\end{eqnarray*}
Consequently, $\nu(A)\le\kappa$ and $\nu(A_1)\le\kappa$. Theorem \ref{thm-1.3}
follows from the Nash embedding theorem as in the analytic category we may take $\kappa\le\frac12m(m+1)$.\hfill\qedbox 

\section*{Acknowledgments} J.C. D{\'\i}az-Ramos and E. Garc{\'\i}a-R{\'\i}o are supported by project
BFM2003-02949, Spain. Research of P. Gilkey partially supported by the
MPI (Leipzig). 

\section*{Dedication}\noindent{\bf 11 de Marzo de 2004 Madrid:} En memoria de todas las
v\'\i ctimas inocentes. Todos \'\i bamos en ese tre{n}. (In memory of all these innocent victims. We were all on that
train.)


\begin{thebibliography}{AAA}
\bibitem{Ni04c} Y. Nikolayevsky, Riemannian manifolds of dimension 7 whose skew-symmetric curvature operator has constant
eigenvalues, math.DG/0311429.
\end{thebibliography}
\end{document}